\documentclass[12pt]{article}
\usepackage{amssymb}
\usepackage{amsmath}

\begin{document}

\begin{center}
\textbf{ON EXTENSION TO FOURIER TRANSFORMS}
\end{center}

\begin{center}
\textsc{Vladimir Lebedev}
\end{center}

\quad

\begin{quotation}
{\small\textsc{Abstract.} For $n$-point sets $K$ in an LCA
group $\Gamma$ we obtain an estimate for the norm of ``the
best'' extension operator from the space $l^\infty(K)$ of
bounded functions on $K$ to the space $A(\Gamma)$ of Fourier
transforms. As a simple consequence our estimate implies that
if $K$ is an infinite closed subset of $\Gamma$ then there does
not exist a bounded linear extension operator from the space
$C_0(K)$ of continuous functions on $K$ vanishing at infinity
to $A(\Gamma)$. The latter result generalizes a result by
Graham who considered the case of compact subsets $K$.


}
\end{quotation}

\quad

\begin{center}
\textsc{1. Introduction and statement of results}
\end{center}

Let $G$ be a locally compact abelian (LCA) group and $\Gamma$
its dual. We consider the space $A(\Gamma)$ of Fourier
transforms, i.e., the space of all functions $f$ of the form
$f=\widehat{\lambda}$, where $\lambda\in L^1(G)$. The norm on
$A(\Gamma)$ is defined by
$\|f\|_{A(\Gamma)}=\|\lambda\|_{L^1(G)}$. Here, $L^1(G)$ is the
space of all complex-valued functions on $G$ integrable with
respect to the Haar measure $m_G$ on $G$ and
$\|\lambda\|_{L^1(G)}=\int_G|\lambda(g)|dm_G(g)$. Note that
$\|\cdot\|_{A(\Gamma)}$ does not depend on the normalization of
$m_G$. Clearly, $A(\Gamma)$ is a Banach space.

For a closed set $K\subseteq\Gamma$ let $C_0(K)$ stand for the
space of all complex-valued continuous functions $f$ on $K$
vanishing at infinity. It is a Banach space with respect to the
norm $\|f\|_{C_0(K)}=\sup_{\gamma\in K}|f(\gamma)|$. We recall
that a function $f$ on $K$ is said to vanish at infinity if for
every $\varepsilon>0$ the set $\{\gamma\in K :
|f(\gamma)|>\varepsilon\}$ is contained in some compact set. If
$K$ is compact, then $C_0(K)$ coincides with the space $C(K)$
of continuous functions on $K$. We invariably use the notation
$C_0(K)$ in both non-compact and compact cases. We have
$A(\Gamma)\subseteq C_0(\Gamma)$ and
$\|f\|_{C_0(\Gamma)}\leq\|f\|_{A(\Gamma)}$.

Closed sets $K\subseteq\Gamma$ with the property that each
function in $C_0(K)$ is a restriction to $K$ of some function
in $A(\Gamma)$, i.e., sets $K$ satisfying
$$
A(\Gamma)_{|K}=C_0(K),
\eqno(1)
$$
are of constant interest in harmonic analysis (we provide
related references and a brief discussion on terminology at the
end of the introduction). It is well-known that each infinite
LCA group has a nontrivial, i.e., infinite, subset $K$
satisfying (1) [5, Theorems 5.6.6 and 5.7.5]. Certainly,
condition (1) just means that for every function $f\in C_0(K)$
there exists an $F\in A(\Gamma)$ which coincides with $f$ on
$K$; one naturally calls $F$ an extension of $f$. At the same
time, Graham [1] showed that if $K$ is an infinite compact
subset of $\Gamma$ then this extension can never be realized by
a bounded linear operator, namely he showed that if $K$ is an
infinite compact subset of $\Gamma$, then there does not exist
a bounded linear operator $\mathcal{E} : C(K)\rightarrow
A(\Gamma)$ such that
$$
(\mathcal{E}f)_{|K}=f
\eqno(2)
$$
for all $f\in C(K)$.

It is worth noting that the proof of this result in [1] uses
tensor algebras.

Generally, given linear spaces $X(K)$ and $Y(\Gamma)$ of
functions on $K$ and $\Gamma$ respectively, we call a linear
operator $\mathcal{E} : X(K)\rightarrow Y(\Gamma)$ an extension
operator if (2) holds for all $f\in X(K)$.

A natural question related to the result of Graham is whether a
similar statement holds in the case of non-compact closed sets
$K$. For instance, let $G$ and $\Gamma$ be the circle group
$\mathbb T$ and the group of integers $\mathbb Z$,
respectively. Consider the subset $K=\{2^j, j=0, 1, 2,
\ldots\}$ of $\mathbb Z$. Being a Sidon set, $K$ has the
property that for each sequence $\{c_j, \,j=0, 1, 2, \ldots\}$
convergent to zero there exists a function $g\in L^1(\mathbb
T)$ with $\widehat{g}(2^j)=c_j, \,j=0, 1, 2, \ldots$. Thus,
condition (1) holds. Does there exist a bounded linear
extension operator from $C_0(K)$ to $A(\mathbb Z)$?

The results of this note are as follows. For a set $K$ let
$l^\infty(K)$ denote the (Banach) space of all complex-valued
bounded functions $f$ on $K$ with
$\|f\|_{l^\infty(K)}=\sup_{\gamma\in K}|f(\gamma)|$. Theorem 1
below is of a quantitative character; it shows, in particular,
that if $K$ is an $n$-point subset of $\Gamma$, then the norm
of any extension operator from $l^\infty(K)$ to $A(\Gamma)$ is
inevitably large as soon as $n$ is large. As a simple and
natural consequence this implies Theorem 2 that establishes
non-existence of a bounded linear extension operator in the
general case of closed sets $K\subseteq\Gamma$, not necessarily
compact ones (thus, in particular, the answer to the above
question about the set $\{2^j, j=0, 1, 2,
\ldots\}\subseteq\mathbb Z$ is negative).

Given a finite set $K\subseteq\Gamma$, define
$$
\alpha_\Gamma(K)\stackrel{\mathrm{def}}{=}\inf\{\|\mathcal{E}\| :
\mathcal{E} \,\, \textrm{is a linear extension operator from} \,\,
l^\infty(K) \,\, \textrm{to} \,\, A(\Gamma)\}.
$$

\quad

\textbf{Theorem 1.} \emph{Let $\Gamma$ be an LCA group. Let $K$
be an $n$-point subset of $\Gamma$. Then
$\sqrt{n/2}\leq\alpha_\Gamma(K)\leq \sqrt{n}$.}

\quad

\textbf{Theorem 2.} \emph{Let $\Gamma$ be an LCA group. Let $K$
be an infinite closed subset of $\Gamma$. Then, there does not
exist a bounded linear extension operator from $C_0(K)$ to
$A(\Gamma)$.}

\quad

For basic results on sets satisfying (1) we refer the reader to
[2--5]. We note that some authors ([2] and [3]) use the term
``Helson set'' for closed sets $K$ with (1), in which case
Sidon sets are just Helson subsets of a discrete group, while
authors of earlier works use the term ``Helson set'' referring
only to compact sets which is the case with the work [1] by
Graham (as well as with the book [5]). To avoid confusion, in
this note we refrain from using any terminology related to
condition (1).

\quad

\begin{center}
\textsc{2. Proofs of the theorems}
\end{center}

We write the Fourier transform $\widehat{\lambda}$ of
$\lambda\in L^1(G)$ as
$$
\widehat{\lambda}(\gamma)=\int_G (-g, \gamma)\lambda(g)dm_G(g),\quad \gamma\in\Gamma,
$$
where $(g, \gamma)$ stands for the value of the character
$\gamma\in\Gamma$ at the point $g\in G$. We actually use $dx$
instead of $dm_G(x)$ in integrals. By $\overline{z}$ we denote
the complex conjugate of a number~$z$.

\quad

\textbf{Proof of Theorem 1.} We begin by proving the lower
bound. The idea of the proof is to consider all functions on
$K$ which take values $\pm 1$ and average the norms of their
images under an extension operator.

Let $r_j, \,j=1, 2, \ldots,$ be the Rademacher functions on the
interval~$[0, 1]$:
$$
r_j(\theta)=\mathrm{sgn}\sin(2^j\pi\theta), \qquad \theta\in[0, 1]
$$
(here $\mathrm{sgn}\,a=a/|a|$ for $a\neq 0$ and
$\mathrm{sgn}\,0=0$). According to the well-known Khinchin
inequality, for an arbitrary $n$ and arbitrary complex numbers
$a_1, a_2, \ldots, a_n$ we have
$$
c\bigg(\sum_{j=1}^n |a_j|^2\bigg)^{1/2}\leq
\int_{[0, 1]}\bigg|\sum_{j=1}^n a_j r_j(\theta)\bigg|d\theta,
\eqno(3)
$$
where $c>0$ does not depend on $n$ and $a_1, a_2, \ldots, a_n$.
It was shown by Szarek [6, Theorem 1, Remark 3] that (3) holds
with $c=1/\sqrt{2}$.

Let $K=\{t_j, \,j=1, 2, \ldots, n\}$, where $t_j$'s are all
distinct. Let $\mathcal{E}$ be an arbitrary linear extension
operator from $l^\infty(K)$ to $A(\Gamma)$.

For each $j=1, 2, \ldots, n$ let $\psi_j$ be a function on $K$
whose value at the point $t_j$ equals $1$ and whose values at
all other points of $K$ equal $0$. Let
$\Psi_j=\mathcal{E}\psi_j$. We have $\Psi_j=\widehat{\Phi_j}$,
where $\Phi_j\in L^1(G)$.

For each $\theta\in[0, 1]$ let $f_\theta$ be the function on
$K$ defined by $f_\theta(t_j)=r_j(\theta), \,j=1, 2, \ldots n$.
We have $\|f_\theta\|_{l^\infty(K)}\leq 1$ and
$$
f_\theta=\sum_{j=1}^n r_j(\theta)\psi_j.
$$
So,
$$
\bigg\|\sum_{j=1}^n r_j(\theta)\Psi_j\bigg\|_{A(\Gamma)}=
\bigg\|\sum_{j=1}^n r_j(\theta)\mathcal{E}\psi_j\bigg\|_{A(\Gamma)}=
\|\mathcal{E}f_\theta\|_{A(\Gamma)}\leq \|\mathcal{E}\|.
$$
This yields
$$
\bigg\|\sum_{j=1}^n r_j(\theta)\Phi_j\bigg\|_{L^1(G)}\leq \|\mathcal{E}\|.
$$
Thus, for all $\theta\in[0, 1]$ we have
$$
\int_G\bigg|\sum_{j=1}^Nr_j(\theta)\Phi_j(x)\bigg|dx\leq \|\mathcal{E}\|.
$$
By integrating this inequality with respect to $\theta\in[0,
1]$ and applying the Khinchin inequality with the Szarek
constant, we obtain
$$
\int_G\frac{1}{\sqrt{2}}\bigg(\sum_{j=1}^n |\Phi_j(x)|^2\bigg)^{1/2}dx\leq \|\mathcal{E}\|.
\eqno(4)
$$

  At the same time for all $x$
$$
\sum_{j=1}^n |\Phi_j(x)|\leq
n^{1/2}\bigg(\sum_{j=1}^n |\Phi_j(x)|^2\bigg)^{1/2},
$$
so (4) implies
$$
\frac{1}{\sqrt{2}}\int_G\sum_{j=1}^n |\Phi_j(x)|dx\leq n^{1/2}\|\mathcal{E}\|.
$$
Thus,
$$
\frac{1}{\sqrt{2}}\sum_{j=1}^n\|\Psi_j\|_{A(\Gamma)}\leq n^{1/2}\|\mathcal{E}\|.
\eqno(5)
$$
It remains to note that, since for each $j$
$$
\|\Psi_j\|_{A(\Gamma)}\geq\sup_{\gamma\in\Gamma}|\Psi_j(\gamma)|\geq
\Psi_j(t_j)=\psi_j(t_j)=1,
$$
estimate (5) yields
$$
\frac{1}{\sqrt{2}}n\leq n^{1/2}\|\mathcal{E}\|,
$$
which completes the proof of the lower bound.

\quad

To prove the upper bound we need the lemma below. This lemma is
a modification of the theorem on local units [5, Theorem
2.6.1].

\quad

\textbf{Lemma.} \emph{Let $V$ be a neighborhood of zero in
$\Gamma$. Then there exists $\Delta\in A(\Gamma)$ such that}

(i)\qquad $\Delta(t)=0$ \emph{for all} $t\notin V$;

(ii)\qquad $\Delta(0)=1$;

(iii)\qquad $0\leq \Delta(t)\leq 1$ \emph{for all}
$t\in\Gamma$;

(iv)\qquad $\Delta=\widehat{\lambda}$, \emph{where} $\lambda\in
L^1(G)$, $\lambda(x)\geq 0$ \emph{for all} $x\in G$,
$\|\lambda\|_{L^1(G)}=1$.

\quad

\textbf{Proof of the Lemma.}  In the ordinary way, assume that
the Fourier transform is extended from $L^1\cap L^2(G)$ to
$L^2(G)$ so that it is a one-to-one mapping of $L^2(G)$ onto
$L^2(\Gamma)$. Assume also that the Haar measure $m_\Gamma$ is
normalized so that
$\|\widehat{\lambda}\|_{L^2(\Gamma)}=\|\lambda\|_{L^2(G)},
\,\lambda\in L^2(G)$. This ensures that
$\widehat{ab}=\widehat{a}\ast \widehat{b}$ for all $a, b\in
L^2(G)$, where $\ast$ stands for convolution (see [5] for
details). Given a set $D\subseteq\Gamma$, by $1_D$ we denote
its indicator function: $1_D(\gamma)=1$ for $\gamma\in D$ and
$1_D(\gamma)=0$ for $\gamma\in \Gamma\setminus D$. One obtains
the functions $\Delta$ and $\lambda$ as follows. Choose a
neighborhood $I$ of zero in $\Gamma$ so that $I-I\subseteq V,
\,0<m_\Gamma(I)<\infty$. Define
$$
\Delta=\frac{1}{m_\Gamma(I)}1_I \ast 1_{-I}.
$$
Clearly, $1_I=\widehat{\xi}$, where $\xi$ is a certain function
in $L^2(G)$, whence $1_{-I}=\widehat{\overline{\xi}}$. Let
$$
\lambda=\frac{1}{m_\Gamma(I)}|\xi|^2.
$$
We have $\lambda\in L^1(G)$ and $\Delta=\widehat{\lambda}$. The
properties (i) -- (iv) are obvious. The lemma is proved.

\quad

To obtain the upper bound in Theorem 1 choose a neighborhood
$V$ of zero in $\Gamma$ so small that
$$
((K-K)\setminus\{0\})\cap V=\varnothing
\eqno(6)
$$
Let $\Delta\in A(\Gamma)$ and $\lambda\in L^1(G)$ be
corresponding functions from the Lemma
($\Delta=\widehat{\lambda}$). Note that (6) implies
$$
\Delta(q-p)=0\quad\textrm{for}\quad p,\,q\in K, \,p\neq q.
\eqno(7)
$$

Define an operator $\mathcal{E} : l^\infty(K)\rightarrow
A(\Gamma)$ as follows. For a function $f$ on $K$ let
$$
\mathcal{E}f(t)=\sum_{\gamma\in K}f(\gamma)\Delta(t-\gamma), \qquad t\in\Gamma.
\eqno(8)
$$
Clearly, $\mathcal{E}$ is a linear extension operator from
$l^\infty(K)$ to $A(\Gamma)$. To estimate its norm consider an
arbitrary function $f$ on $K$ with $\|f\|_{l^\infty(K)}=1$. We have
$\mathcal{E}f=\widehat{S}$ where
$$
S(x)=\sum_{\gamma\in K}f(\gamma)(x, \gamma)\lambda(x), \qquad x\in G.
$$
So,
$$
\|\mathcal{E}f\|_{A(\Gamma)}=\|S\|_{L^1(G)}=
\int_G\bigg|\sum_{\gamma\in K}f(\gamma)(x, \gamma)\bigg|\lambda(x)dx.
\eqno(9)
$$

Note that the measure $m$ given by $dm(x)=\lambda(x)dx$ is a
probability measure on $G$, hence,
$$
\int_G\bigg|\sum_{\gamma\in K}f(\gamma)(x, \gamma)\bigg|\lambda(x)dx\leq
\bigg(\int_G\bigg|\sum_{\gamma\in K}f(\gamma)(x, \gamma)\bigg|^2\lambda(x)dx\bigg)^{1/2}.
$$
Thus (see (9), (7))
\begin{equation*}
\begin{split}
\|\mathcal{E}f\|_{A(\Gamma)}^2 & \leq \int_G\bigg|\sum_{\gamma\in K}f(\gamma)(x, \gamma)\bigg|^2\lambda(x)dx \\
& = \int_G\bigg(\sum_{p, q\in K}f(p)\overline{f(q)}(x, p-q)\bigg)\lambda(x)dx \\
& = \sum_{p, q\in K}f(p)\overline{f(q)}\bigg(\int_G(x,
p-q)\lambda(x)dx\bigg) \\
& = \sum_{p, q\in
K}f(p)\overline{f(q)}\Delta(q-p)=\sum_{\gamma\in
K}|f(\gamma)|^2\leq n.
\end{split}
\end{equation*}
This completes the proof of Theorem 1.

\quad

\textbf{Proof of Theorem 2.} First note that if
$K\subseteq\Gamma$ is a finite set, then there exists a linear
extension operator from $l^\infty(K)$ to $C_0(\Gamma)$ whose
norm is equal to $1$. Indeed, the operator defined by (8) has
this property. Let now $K$ be an infinite closed subset of
$\Gamma$. For each $n=1, 2, \ldots$ chose an $n$-point subset
$K_n$ of $K$. Let $\Sigma_n : l^\infty(K_n)\rightarrow
C_0(\Gamma)$ be a linear extension operator with norm equal to
$1$. Let $\mathcal{R} : C_0(\Gamma)\rightarrow C_0(K)$ stand
for the operator of restriction to $K$, i.e., the operator that
takes $f\in C_0(\Gamma)$ to its restriction $f_{|_K}$. Assuming
that there exists a bounded linear extension operator
$\mathcal{E} : C_0(K)\rightarrow A(\Gamma)$, consider the
operator $\mathcal{E}\mathcal{R}\Sigma_n$. This operator is a
linear extension operator from $l^\infty(K_n)$ to $A(\Gamma)$
and its norm is at most $\|\mathcal{E}\|$. This contradicts
Theorem 1 if $n$ is large enough.

\quad

\textbf{Remark.} Let $M(G)$ be the space of all complex regular
bounded measures on $G$. Consider the space $B(\Gamma)$ of
Fourier--Stieltjes transforms, i.e, the space of all functions
$f$ of the form $f=\widehat{\mu}$, where $\mu\in M(G)$ ([5,
Section 1.3]). The norm on $B(\Gamma)$ is defined by
$\|f\|_{B(\Gamma)}=\|\mu\|_{M(G)}$. Clearly, $B(\Gamma)$ is a
Banach space and $A(\Gamma)\subseteq B(\Gamma)$  (these spaces
coincide when $\Gamma$ is compact). Theorem 2 can be
strengthened by replacing $A(\Gamma)$ with $B(\Gamma)$. To see
this, we define for a finite set $K\subseteq\Gamma$
$$
\beta_\Gamma(K)\stackrel{\mathrm{def}}{=}\inf\{\|\mathcal{E}\| :
\mathcal{E} \,\, \textrm{is a linear extension operator from} \,\,
l^\infty(K) \,\, \textrm{to} \,\, B(\Gamma)\}.
$$
Note that
$$
\alpha_\Gamma(K)=\beta_\Gamma(K).
\eqno(10)
$$
Indeed, relation $\alpha_\Gamma(K)\geq\beta_\Gamma(K)$ is
obvious. To verify that $\alpha_\Gamma(K)\leq\beta_\Gamma(K)$,
we let $\delta>0$ and consider a function $\chi_\delta\in
A(\Gamma)$ such that $\chi_\delta=1$ on $K$ and
$\|\chi_\delta\|_{A(\Gamma)}\leq 1+\delta$ (see [3, Proposition
A.5.1]). Assuming that $\mathcal{E} : C_0(K)\rightarrow
B(\Gamma)$ is an extension operator, we let $\mathcal{E}_\delta
f=\chi_\delta\mathcal{E}f$ for $f\in C_0(K)$. Note that if
$a\in A(\Gamma), b\in B(\Gamma)$, then $ab\in A(\Gamma)$ and
$\|ab\|_{A(\Gamma)}\leq\|a\|_{A(\Gamma)}\|b\|_{B(\Gamma)}$. So,
$\mathcal{E}_\delta$ is an extension operator from $C_0(K)$ to
$A(\Gamma)$ and $\|\mathcal{E}_\delta\|_{C_0(K)\rightarrow
A(\Gamma)}\leq (1+\delta)\|\mathcal{E}\|_{C_0(K)\rightarrow
B(\Gamma)}$. Thus, we obtain (10). It remains to proceed as in
the proof of Theorem 2 with obvious modifications.

\quad

\quad

\quad

\quad

\quad

\begin{center}
\textsc{References}
\end{center}

\flushleft
\begin{enumerate}

\item C. C. Graham, ``Helson sets and simultaneous extension
    to Fourier transforms'', \emph{Studia Math.} 43 (1972),
    57--59.

\item C. C. Graham, K. E. Hare, \emph{Interpolation and Sidon
    Sets for Compact Groups}, CMS Books in Mathematics,
    Springer, Boston, MA, 2013.

\item  C. C. Graham, O. C. McGehee, \emph{Essays in
    commutative harmonic analysis}, Grundlehren der
    mathematischen Wissenschaften, vol. 238, Springer,
    Berlin, 1979.

\item J.-P. Kahane, \emph{S\'erie de Fourier absolument
    convergentes}, Springer-Verlag, Berlin--Heidelberg--New
    York, 1970.

\item W. Rudin, \emph{Fourier analysis on groups},
    Interscience Publishers, New York--London, 1962.

\item  S. J. Szarek, ``On the best constants in the Khinchin
    inequality'', Studia Math. 58 (1976), 197--208.

\end{enumerate}

\quad

\qquad School of Applied Mathematics\\
\qquad Natioinal Research University Higher School of Economics\\
\qquad (HSE University)\\
\qquad 34 Tallinskaya St.\\
\qquad Moscow, 123458 Russia\\

\quad

\qquad e-mail: \emph{lebedevhome@gmail.com}

\end{document}